\documentclass[fleqn]{mat01}
\usepackage{times,mathtimy,amssymb,latexsym}
\begin{document}

\setcounter{page}{349} \firstpage{349}

\newtheorem{theo}{\bf Theorem}
\newtheorem{coro}[defin]{\rm COROLLARY}
\newtheorem{lem}{Lemma}
\newtheorem{pol}{\it Proof of Lemma}
\newtheorem{pot}{\it Proof of Theorem}

\def\remar{\trivlist \item[\hskip \labelsep{\it Remark.}]}

\def\d{\mbox{\rm d}}
\def\e{\mbox{\rm e}}

\title{On an inequality concerning the polar derivative\\ of a polynomial}

\markboth{A Aziz and N A Rather}{Polar derivative of a polynomial}

\author{A AZIZ and N A RATHER$^{*}$}

\address{Post Graduate Department of Mathematics,
University of Kashmir, Hazratbal, Srinagar~190~006, India\\
\noindent $^{*}$Post Graduate Department of Mathematics,
University of Kashmir, Hazratbal, Srinagar~190~006,
India\\
\noindent E-mail: nisararather@yahoo.co.in}

\volume{117}

\mon{August}

\parts{3}

\pubyear{2007}

\Date{MS received 16 March 2006; revised 24 October 2006}

\begin{abstract}
In this paper, we present a correct proof of an $L_p$-inequality
concerning the polar derivative of a polynomial with restricted
zeros. We also extend Zygmund's inequality to the polar derivative
of a polynomial.
\end{abstract}


\maketitle

\section{Introduction and statement of results}

Let $P(z)$ be a polynomial of degree$ n$ and let ${P}'(z)$ be its
derivative. Then according to the famous result known as
Bernstein's inequality (see \cite{7} or \cite{10})
\begin{equation}
\max_{\vert z\vert = 1} |{{P}'(z)}| \le n\max_{\vert z\vert = 1}
\left| {P(z)} \right|.
\end{equation}
Inequality (1) is sharp and equality in (1) holds for$P(z) =
az^n,a \ne 0$. Inequality (1) was extended to $L_p$-norm by
Zygmund \cite{11} who proved that if $P(z)$ is a polynomial of
degree $n$, then for $p \ge 1$,
\begin{equation}
\left\{ {\int_0^{2\pi } {| {{P}'(\e^{i\theta })}|^p\d\theta}}
\right\}^{1/p} \le n\left\{ {\int_0^{2\pi } {| {P(\e^{i\theta
})}|^p\d\theta } } \right\}^{1/p}.
\end{equation}
The result is sharp and equality in (2) holds for $P(z) = az^n,a
\ne 0$. If we let $p \to \infty$ in (2), we get inequality (1).

Let $D_\alpha P(z)$ denote the polar differentiation of polynomial
$P(z)$ with respect to a real or complex number $\alpha$. Then
\begin{align*}
D_\alpha P(z) = nP(z) + (\alpha - z){P}'(z).
\end{align*}
The polynomial $D_\alpha P(z)$is of degree at most $n-1$ and it
generalizes the ordinary derivative in the sense that
\begin{equation*}
\mathop {\lim}\limits_{\alpha \to \infty } \frac{D_\alpha
P(z)}{\alpha } = {P}'(z)
\end{equation*}
uniformly on compact subsets of $C$.\newpage

As an extension of (1) to the polar derivative, Aziz and Shah
(Theorem~4 with $k = 1$, \cite{3}) have shown that if $P(z)$ is a
polynomial of degree $n$, then for every complex number $\alpha$
with $\vert \alpha \vert \ge 1$,
\begin{equation}
|{D_\alpha P(z)}| \le n\vert \alpha \vert \max_{\vert z\vert = 1}
|{P(z)}|\quad\hbox{for}\ \vert z\vert  = 1.
\end{equation}
Inequality (3) becomes equality for $P(z) = az^n, a \ne 0$.

If we divide the two sides of (3) by $|\alpha|$ and let $|\alpha|
\to \infty$, we get inequality (1).

It is natural to seek $L_p$-norm analog of inequality (3). In view
of the $L_p$-norm extension (2) of inequality (1), one would
expect that if $P(z)$ is a polynomial of degree $n$, then
\begin{equation}
\left\{ {\int_0^{2\pi } {| {D_\alpha P(\e^{i\theta })} |^p
\d\theta } } \right\}^{1/p} \le n| \alpha |\left\{ {\int_0^{2\pi }
{| {P(\e^{i\theta })} |^p }}\right\}^{1/p}
\end{equation}
will be $L_p$-norm extension of (3) analogous to (2). But
unfortunately inequality (4) is not, in general, true for every
real or complex number $\alpha$. To see this, we take in
particular $p = 2, P(z) = (1 - iz)^n$ and $\alpha  = i\beta $
where $\beta$ is any positive real number such that
\begin{equation}
1 \le \beta < \frac{n + \sqrt {2n(2n - 1)} }{3n - 2}.
\end{equation}
Now
\begin{align*}
D_\alpha P(z) &= n(1 - iz)^n - ni(\alpha - z)(1 - iz)^{n - 1}\\[.2pc]
 &= n(1 - iz)^{n - 1}(1 - i\alpha )
\end{align*}
so that
\begin{align}
\int_0^{2\pi } \vert D_\alpha P(\e^{i\theta })\vert^p \d\theta  &=
n^2\vert 1 - i\alpha \vert ^2\int_0^{2\pi }
{\vert 1 - i\e^{i\theta }\vert^{2(n - 1)}\d\theta }\nonumber\\[.4pc]
&= n^2\vert 1 - i\alpha \vert^2\int_0^{2\pi } {\vert (1 -
i\e^{i\theta })^{n - 1}\vert^2 \d\theta }\nonumber\\[.4pc]
&= n^2\vert 1 - i\alpha \vert^2\int_0^{2\pi }
\left\vert\begin{pmatrix} n - 1\\[.1pc] 0 \end{pmatrix} -
\begin{pmatrix} n - 1\\[.1pc] 1 \end{pmatrix}
(i\e^{i\theta })\right.\nonumber\\[.4pc]
&\quad\, + \begin{pmatrix} n - 1\\[.1pc] 2
\end{pmatrix} (i\e^{i\theta })^2 - \cdots\nonumber\\[.4pc]
&\quad\, \left.+ ( - 1)^{n - 1}\begin{pmatrix} n - 1\\[.1pc] n-1
\end{pmatrix} (i\e^{i\theta })^{n - 1}\right\vert ^2 \d\theta\nonumber\\[.4pc]
&= 2\pi n^2\vert 1 - i\alpha \vert ^2\left(
\begin{pmatrix} n - 1\\[.1pc] 0 \end{pmatrix}^{2} +
\begin{pmatrix} n - 1\\[.1pc] 1 \end{pmatrix}^{2} \right.\nonumber\\[.4pc]
&\quad\, \left.+\begin{pmatrix} n - 1\\[.1pc] 2 \end{pmatrix}^{2} +
\cdots +
\begin{pmatrix} n - 1\\[.1pc] n-1 \end{pmatrix}^{2}
\right)\nonumber\\[.4pc]
&= 2\pi n^2\vert 1 - i\alpha \vert^2
\begin{pmatrix} 2(n-1)\\[.1pc] n-1 \end{pmatrix}.
\end{align}
Also,
\begin{align}
n^2\vert \alpha \vert^2\int_0^{2\pi } \vert P(\e^{i\theta
})\vert^p \d\theta &= n^2\vert \alpha \vert^2\int_0^{2\pi }
{\vert 1 - i\e^{i\theta }\vert ^{2n}\d\theta } \nonumber\\[.4pc]
&= n^2\vert \alpha \vert ^2\int_0^{2\pi } {\vert (1 -
i\e^{i\theta })^n\vert^2 \d\theta }\nonumber\\[.4pc]
 &= n^2\vert \alpha \vert^2\int_0^{2\pi }
\begin{pmatrix} n\\[.1pc] 0 \end{pmatrix} - \begin{pmatrix} n\\[.1pc]
1 \end{pmatrix}
(i\e^{i\theta }) + \begin{pmatrix} n\\[.1pc] 2 \end{pmatrix}
(i\e^{i\theta })^2 \nonumber\\[.4pc]
&\quad\, - \cdots + ( - 1)^n \begin{pmatrix} n\\[.1pc] n \end{pmatrix}
(i\e^{i\theta })^n\vert ^2 \d\theta\nonumber\\[.4pc]
 &= 2\pi n^2\vert \alpha \vert ^2\left(
\begin{pmatrix} n\\[.1pc] 0 \end{pmatrix}^{2} +
\begin{pmatrix} n\\[.1pc] 1 \end{pmatrix}^{2} + \cdots +
\begin{pmatrix} n\\[.1pc] n \end{pmatrix}^{2}\right)\nonumber\\[.4pc]
&= 2\pi n^2\vert \alpha \vert^2\begin{pmatrix} 2n\\[.1pc] n
\end{pmatrix}.
\end{align}
Using (6) and (7) in (4), we get
\begin{equation*}
2\pi n^2\begin{pmatrix} 2(n-1)\\[.1pc] n-1 \end{pmatrix}\vert 1 - i\alpha
\vert^2 \le 2\pi n^2\vert \alpha \vert ^2\begin{pmatrix} 2n \\[.1pc]
n
\end{pmatrix}.
\end{equation*}
This implies
\begin{equation}
n\vert 1 - i\alpha \vert^2 \le 2(2n - 1)\vert \alpha \vert ^2.
\end{equation}
Setting $\alpha = i\beta$ in (8), we get
\begin{equation*}
n(1 + \beta )^2 \le 2(2n - 1)\beta ^2.
\end{equation*}
This inequality can be written as
\begin{equation}
\left( {\beta - \frac{n + \sqrt {2n(2n - 1)} }{3n - 2}}
\right)\left( {\beta - \frac{n - \sqrt {2n(2n - 1)} }{3n - 2}}
\right) \ge 0.
\end{equation}
Since $\beta \ge 1$, we have
\begin{align*}
\left( {\beta - \frac{n - \sqrt {2n(2n - 1)} }{3n - 2}} \right)
&\ge 1 - \frac{n - \sqrt {2n(2n - 1)} }{3n - 2}\\[.4pc]
 &= \frac{2(n - 1) + \sqrt {2n(2n - 1)} }{3n - 2} > 0
\end{align*}
and hence from (9), it follows that
\begin{equation*}
\left({\beta - \frac{n + \sqrt {2n(2n - 1)} }{3n - 2}} \right) \ge
0.
\end{equation*}
This gives
\begin{equation*}
\beta \ge \frac{n + \sqrt {2n(2n - 1)} }{3n - 2},
\end{equation*}\pagebreak

\noindent which clearly contradicts (5). Hence inequality (4) is
not, in general, true for all polynomials $P(z)$ of degree $n \ge
1$.

However, we have been able to prove the following generalization of (2) to
the polar derivatives.

\begin{theo}[\!]
If $P(z)$ is a polynomial of degree $n${\rm ,} then for every
complex number $\alpha$ and $p \ge 1${\rm ,}
\begin{equation}
\left\{ {\int_0^{2\pi } {| {D_\alpha P(\e^{i\theta })} |^p
\d\theta } } \right\}^{1/p} \le n(| \alpha | + 1)\left\{
{\int_0^{2\pi } {| {P(\e^{i\theta })} |^p \d\theta } }
\right\}^{1/p}.
\end{equation}\vspace{-.3pc}
\end{theo}

\begin{remar}
If we divide both sides of (10) by $\vert\alpha\vert $ and make
$\vert \alpha \vert  \to \infty $, we get inequality (2) due to
Zygmund \cite{11}.

For polynomials $P(z)$ which does not vanish in the unit disk, the
right-hand side of (2) can be improved. In fact , in this
direction, it was shown by De-Bruijn \cite{4} that if $P(z)$ does
not vanish in $\vert z\vert < 1$, then for $p \ge 1$,
\begin{equation}
\left\{ {\int_0^{2\pi } {| {{P}'(\e^{i\theta })} |^p \d\theta } }
\right\}^{1/p} \le nC_p \left\{ {\int_0^{2\pi } {|{P(\e^{i\theta
})} |^p \d\theta } } \right\}^{1/p},
\end{equation}
where
\begin{equation}
C_p = \left\{ {\frac{1}{2\pi }\int_0^{2\pi } {| {1 + \e^{i\beta }}
|^p \d\beta } } \right\}^{ - 1/p}.
\end{equation}
Inequality (11) is best possible with equality for $P(z) = az^n +
b,\vert a\vert = \vert b\vert $. If we let $p  \to \infty$ in
(11), it follows that if $P(z) \ne  0$ for $\vert z\vert <1$, then
\begin{equation}
\max_{\vert z\vert = 1} | {{P}'(z)} | \le \frac{n}{2} \max_{\vert
z\vert = 1} | {P(z)} |.
\end{equation}
Inequality (13) was conjectured by Erd\"{o}s and later verified by
Lax \cite{6}. Aziz \cite{1} extended (13) to the polar derivative
of a polynomial and proved that if $P(z)$ is a polynomial of
degree $n$ which does not vanish in $\vert z\vert < 1$, then for
every complex number $\alpha$ with $\vert \alpha \vert \ge 1$,
\begin{equation}
\max_{\vert z\vert = 1} | {D_\alpha P(z)} | \le \frac{n}{2}(\vert
\alpha \vert + 1)\max_{\vert z\vert = 1} | {P(z)} |.
\end{equation}
The estimate (14) is best possible with equality for $P(z) = z^n +
1$. If we divide both sides of (14) by $\vert \alpha\vert$ and
make $\vert \alpha \vert \to \infty$, we get inequality (13) due
to Lax \cite{6}.

While seeking the desired extension to the polar derivatives,
recently Govil {\it et~al} \cite{5} have made an incomplete
attempt by claiming to have proved the following generalization of
(11) and (14).
\end{remar}

\begin{theo}[\!]
If $P(z)$ is a polynomial of degree n which does not vanish in
\hbox{$\vert z\vert < 1$}{\rm ,} then for every complex number
$\alpha $ with $\vert  \alpha \vert  \ge  1$ and $p \ge 1${\rm ,}
\begin{equation}
\left\{ {\int_0^{2\pi } {| {D_\alpha P(\e^{i\theta })} |^p
\d\theta } } \right\}^{1/p} \le n(| \alpha | + 1)C_p \left\{
{\int_0^{2\pi } {|{P(\e^{i\theta })}|^p \d\theta } }
\right\}^{1/p},
\end{equation}
where $C_p $ is defined by $(12)$.
\end{theo}

Unfortunately the proof of this theorem, which is the main result
(Theorem~1.1 of \cite{5}) given by Govil, Nyuydinkong and Tameru
is not correct, because the claim made by the authors on page 624
in lines 12 to 16 by using Lemma~2.3 is incorrect. The reason
being that their polynomial
\begin{equation*}
D_\alpha P_n (z) + \e^{i\gamma }\{ {n\bar{\alpha }zP_n (z) + (1 -
\bar{\alpha }z)z{P}'_n (z)} \}, \quad z = \e^{i\theta},
\end{equation*}
in general does not take the form
\begin{equation*}
\sum\limits_{k = 0}^n {l_k a_k z^k}, \quad z = \e^{i\theta }
\end{equation*}
where
\begin{align*}
P_n (z) = \sum\limits_{k = 0}^n {a_k z^k}
\end{align*}
and the complex numbers $l_k $ defined by them on page 624, line
10, by
\begin{equation*}
L( {P_n (\e^{i\theta })}) = [ {\Lambda P_n (\e^{i\theta })}
]_{\theta = 0} = \sum\limits_{k = 0}^n {l_k a_k }
\end{equation*}
along with the equation (24) of \cite{5}.

It is worthwhile to note here that if we take
\begin{equation*}
L( {P_n (\e^{i\theta })} ) = [\,nP_n (\e^{i\theta }) + (\alpha -
\e^{i\theta }){P}'_n (\e^{i\theta })]_{\theta = 0}
\end{equation*}
and use the same argument as used by Govil {\it et~al} (page~624,
line 10 of \cite{5}), then in view of the inequality
\begin{equation*}
| {D_\alpha P(z)} | \le n\vert \alpha \vert \max_{\vert z\vert =
1} | {P(z)} |\quad\hbox{for}\ \vert z\vert  = 1
\end{equation*}
(see Theorem~4 with $k = 1$ of \cite{3}), the above bounded
functional has norm $N \le n\vert \alpha \vert $. Therefore, if we
use Lemma~2.3 of \cite{5} which is due to Rahman (Lemma~3 of
\cite{8}), it would follow that
\begin{equation*}
\left\{ {\int_0^{2\pi } {| {D_\alpha P(\e^{i\theta })} |^p
\d\theta } } \right\}^{1/p} \le n| \alpha | \left\{ {\int_0^{2\pi
} {| {P(\e^{i\theta })} |^p \d\theta } } \right\}^{1/p}
\end{equation*}
for every $p \ge  1$ and $\vert \alpha\vert \ge 1$, which is not
true in general as shown above.

Here we shall also present a correct proof of Theorem~2, which
shall validate Theorems~1.2 and 1.3 of Govil {\it et~al} \cite{5}
as well. Finally we shall also present a short proof of
Theorem~1.3 of \cite{5}. That is, we prove the following.

\begin{theo}[\!]
If $P(z)$ is a self-inversive polynomial of degree $n${\rm ,} then
for every complex number $\alpha $ and $p \ge 1${\rm ,}
\begin{equation}
\left\{ {\int_0^{2\pi } {| {D_\alpha P(\e^{i\theta })}|^p \d\theta
} } \right\}^{1/p} \le n(| \alpha | + 1)C_p \left\{ {\int_0^{2\pi
} {| {P(\e^{i\theta })} |^p \d\theta } } \right\}^{1/p},
\end{equation}
where $C_p$ is the same as in Theorem~$2$.
\end{theo}

\section{Lemmas}

For the proofs of these theorems, we need the following
lemmas.\pagebreak

\setcounter{lem}{0}
\begin{lem}
If $P(z)$ is a polynomial of degree $n$ which does not vanish in
$\vert z\vert  < 1${\rm ,} and $Q(z) = z^n\overline{P(1/\bar
{z})}${\rm ,} then for every complex number $\alpha$ with $\vert
\alpha \vert \ge 1${\rm ,}
\begin{equation*}
| {D_\alpha P(z)} | \le | {D_\alpha Q(z)} |\quad\hbox{for}\ \vert
z\vert \ge  1.
\end{equation*}
\end{lem}
Lemma~1 is due to Aziz (p.~190 of \cite{1}).

\begin{lem}
If $P(z)$ is a polynomial of degree $n$ and $Q(z) = z^n\overline
{P(1/\bar {z})}${\rm ,} then for every $p \ge 0$ and $\beta$
real{\rm ,}
\begin{equation*}
\int_0^{2\pi } \int_0^{2\pi } {| {{Q}'(\e^{i\theta }) + \e^{i\beta
}{P}'(\e^{i\theta })} |}^p \d\theta \d\beta \le 2\pi
n^p\int_0^{2\pi } {|{P(\e^{i\theta })} |^p \d\theta }.
\end{equation*}
\end{lem}
Lemma~2 is due to Aziz \cite{2} (see also \cite{8}). We also need
the following lemma.

\begin{lem}
If $P(z)$ is a polynomial of degree $n, P(0) \neq 0$ and $Q(z) =
z^n\overline{P(1/\bar {z})}${\rm ,} then for every complex number
$\alpha${\rm ,} $p \ge 1$ and $\beta$ real{\rm ,}
\begin{equation*}
\hskip -4pc \int_0^{2\pi } \int_0^{2\pi } {| {D_\alpha
Q(\e^{i\theta }) + \e^{i\beta }D_\alpha P(\e^{i\theta })} |} ^p
\d\theta \d\beta \le 2\pi n^p(\vert \alpha \vert +
1)^p\int_0^{2\pi } {| {P(\e^{i\theta })}|^p \d\theta }.
\end{equation*}
\end{lem}

\setcounter{pol}{2}
\begin{pol}
{\rm We have by Minkowski's inequality for every $p \ge 1$ and
$\beta$ real,
\begin{align}
&\left\{ {\int_0^{2\pi } {\int_0^{2\pi } {\vert D_\alpha
Q(\e^{i\theta }) + \e^{i\beta }D_\alpha P(\e^{i\theta })} \vert ^p
\d\theta \d\beta } } \right\}^{1/p}\nonumber\\[.4pc]
&\quad\, = \left\lbrace \int_0^{2\pi } \int_0^{2\pi } {\vert (
{nQ(\e^{i\theta }) + (\alpha - \e^{i\theta }){Q}'(\e^{i\theta })}
)} + \e^{i\beta }( nP(\e^{i\theta })\right.\nonumber\\[.4pc]
&\qquad\,\left.\phantom{\int_0^{2\pi }}\hskip -1.87pc + (\alpha -
\e^{i\theta }){P}'(\e^{i\theta })
)\vert ^p \d\theta \d\beta\right\rbrace^{1/p}\nonumber\\[.4pc]
&\quad\,= \left\lbrace \int_0^{2\pi } \int_0^{2\pi } {\vert (
{nQ(\e^{i\theta }) - \e^{i\theta }{Q}'(\e^{i\theta })} )} +
\e^{i\beta }( {nP(\e^{i\theta }) - \e^{i\theta }{P}'(\e^{i\theta
})} )\right.\nonumber\\[.4pc]
&\qquad\,\left.\phantom{\int_0^{2\pi }}\hskip -1.87pc + \alpha (
{Q}'(\e^{i\theta }) + \e^{i\beta }{P}'(\e^{i\theta }) )\vert ^p
\d\theta \d\beta
\right\rbrace^{1/p}\nonumber\\[.4pc]
&\quad\, \le \left\lbrace \int_0^{2\pi } \int_0^{2\pi } {\vert (
{nQ(\e^{i\theta }) - \e^{i\theta }{Q}'(\e^{i\theta })} )}\right.\nonumber\\[.4pc]
&\qquad\,\left.\phantom{\int_0^{2\pi }}\hskip -1.87pc  +
\e^{i\beta }( {nP(\e^{i\theta }) - \e^{i\theta }{P}'(\e^{i\theta
})} )\vert ^p \d\theta \d\beta \right\rbrace^{1/p}\nonumber\\[.4pc]
&\qquad\, + \vert \alpha \vert \left\{ {\int_0^{2\pi }
\int_0^{2\pi } {| {{Q}'(\e^{i\theta }) + \e^{i\beta
}{P}'(\e^{i\theta })} |} ^p \d\theta \d\beta } \right\}^{1/p}.
\end{align}
Since $Q(z) = z^n\overline{P(1/\bar {z})}$, we have $P(z) =
z^n\overline{Q(1/\bar {z})}$ and it can be easily verified that
for $0  \le  \theta  < 2\pi$,
\begin{equation}
nP(\e^{i\theta }) - \e^{i\theta }{P}'(\e^{i\theta }) = \e^{i(n -
1)\theta }\overline {{Q}'(\e^{i\theta }} )
\end{equation}
and
\begin{equation}
nQ(\e^{i\theta }) - \e^{i\theta }{Q}'(\e^{i\theta }) = \e^{i(n -
1)\theta }\overline {{P}'(\e^{i\theta })}.
\end{equation}
Using (18) and (19) in (17), we obtain
\begin{align*}
&\left\{ {\int_0^{2\pi } {\int_0^{2\pi } {\vert D_\alpha
Q(\e^{i\theta }) + \e^{i\beta }D_\alpha P(\e^{i\theta })} \vert ^p
\d\theta \d\beta } } \right\}^{1/p}\\[.4pc]
&\quad\, \le \left\{ {\int_0^{2\pi } \int_0^{2\pi } {| {\e^{i(n -
1)\theta }\overline {{P}'(\e^{i\theta })} + \e^{i\beta }\e^{i(n -
1)\theta }\overline {{Q}'(\e^{i\theta })} } |} ^p \d\theta \d\beta
} \right\}^{1/p}\\[.4pc]
&\qquad\, + \vert \alpha \vert \left\{ {\int_0^{2\pi }
\int_0^{2\pi } {| {{Q}'(\e^{i\theta }) + \e^{i\beta
}{P}'(\e^{i\theta })} |} ^p \d\theta \d\beta } \right\}^{1/p}\\[.4pc]
&\quad\, = \left\{ {\int_0^{2\pi } \int_0^{2\pi } {|
{{Q}'(\e^{i\theta }) + \e^{i\beta }{P}'(\e^{i\theta })} |}^p
\d\theta \d\beta } \right\}^{1/p}\\[.4pc]
&\qquad\, + \vert \alpha \vert \left\{ {\int_0^{2\pi }
\int_0^{2\pi } {| {{Q}'(\e^{i\theta }) + \e^{i\beta
}{P}'(\e^{i\theta })} |} ^p \d\theta \d\beta } \right\}^{1/p}\\[.4pc]
&\quad\, = (\vert \alpha \vert + 1)\left\{ {\int_0^{2\pi }
\int_0^{2\pi } {| {{Q}'(\e^{i\theta }) + \e^{i\beta
}{P}'(\e^{i\theta })} |} ^p \d\theta \d\beta } \right\}^{1/p}.
\end{align*}
This gives with the help of Lemma~2,
\begin{equation*}
\hskip -4pc\int_0^{2\pi } \int_0^{2\pi } {| {D_\alpha
Q(\e^{i\theta }) + \e^{i\beta }D_\alpha P(\e^{i\theta })} |} ^p
\d\theta \d\beta \le 2\pi n^p(\vert \alpha \vert +
1)^p\int_0^{2\pi } {| {P(\e^{i\theta })} |^p \d\theta} .
\end{equation*}
This completes the proof of Lemma~3.}
\end{pol}

\section{Proofs of the theorems}

\setcounter{pot}{0}
\begin{pot}
{\rm By Lemma~3, we have for every complex number $\alpha$, $p \ge
1$ and $\beta$ real,
\begin{equation}
\hskip -4pc\int_0^{2\pi } \int_0^{2\pi } {| {D_\alpha
Q(\e^{i\theta }) + \e^{i\beta }D_\alpha P(\e^{i\theta })} |} ^p
\d\theta \d\beta \le 2\pi n^p(\vert \alpha \vert +
1)^p\int_0^{2\pi } {| {P(\e^{i\theta })} |^p \d\theta }.
\end{equation}
Using in (20) the fact that for any $p \ge  0$,
\begin{equation*}
\int_0^{2\pi } {\vert a + b\e^{i\beta }\vert ^p \d\beta  \ge 2\pi
\max\{ {\vert a\vert ^p,\vert b\vert ^p} \}}
\end{equation*}
(see inequality (19) of \cite{4}), we obtain
\begin{equation*}
\left\{ {\int_0^{2\pi } {| {D_\alpha P(\e^{i\theta })} |} ^p
\d\theta } \right\}^{1/p} \le n(\vert \alpha \vert + 1)\left\{
{\int_0^{2\pi } {| {P(\e^{i\theta })} |^p \d\theta } }
\right\}^{1/p}.
\end{equation*}
This completes the proof of Theorem~1.}
\end{pot}

\begin{pot}
{\rm Since $P(z)$ is a polynomial of degree $n$ which does not
vanish in $\vert z\vert < 1$, by Lemma~1 we have for every complex
number $\alpha$ with $\vert \alpha \vert \ge 1$,
\begin{equation}
|{D_\alpha P(z)}| \le |{D_\alpha Q(z)}|,\quad\hbox{for}\ \vert
z\vert =1
\end{equation}
where $Q(z) = z^n\overline {P(1/\bar {z})}$. Also by Lemma~3 for
every complex number $\alpha, p \ge 1$ and $\beta$ real,
\begin{equation}
\hskip -4pc \int_0^{2\pi } \left\{ {\int_0^{2\pi } {| {D_\alpha
Q(\e^{i\theta }) + \e^{i\beta }D_\alpha P(\e^{i\theta })} |}^p
\d\beta } \right\}\d\theta  \le 2\pi n^p(\vert \alpha \vert +
1)^p\int_0^{2\pi } {| {P(\e^{i\theta })} |^p \d\theta }.
\end{equation}
Now for every real $\beta$ and $r \ge 1$, we have
\begin{equation*}
\vert r + \e^{i\beta }\vert  \ge \vert 1 + \e^{i\beta }\vert,
\end{equation*}
which implies
\begin{equation*}
\int_0^{2\pi } \vert r + \e^{i\beta }\vert ^p \d\beta  \ge
\int_0^{2\pi } \vert 1 + \e^{i\beta }\vert ^p \d\beta,\quad p \ge
0.
\end{equation*}
If $D_\alpha P(\e^{i\theta }) \ne 0$, we take $r = \vert D_\alpha
Q(\e^{i\theta })\vert/\vert D_\alpha P(\e^{i\theta })\vert$, and
by (21) $r  \ge  1$ and we get
\begin{align*}
&\int_0^{2\pi } {| {D_\alpha Q(\e^{i\theta }) + \e^{i\beta
}D_\alpha P(\e^{i\theta })} |} ^p \d\beta\\[.4pc]
&\quad\,= \vert D_\alpha P(\e^{i\theta })\vert ^p\int_0^{2\pi }
{\left| {\frac{D_\alpha Q(\e^{i\theta })}{D_\alpha P(\e^{i\theta
})} + \e^{i\beta }}
\right|^p \d\beta }\\[.4pc]
 &\quad\,= \vert D_\alpha P(\e^{i\theta })\vert ^p\int_0^{2\pi } {\left|
{\left| {\frac{D_\alpha Q(\e^{i\theta })}{D_\alpha P(\e^{i\theta
})}} \right| + \e^{i\beta }} \right|^p \d\beta }\\[.4pc]
&\quad\,\ge \vert D_\alpha P(\e^{i\theta })\vert^p \int_0^{2\pi }
{| {1 + \e^{i\beta }} |^p \d\beta }.
\end{align*}
For $D_\alpha P(\e^{i\theta }) = 0$, this inequality is trivially
true. Using this in (22), we conclude that for every complex
number $\alpha$ with $\vert \alpha \vert \ge 1$ and $p \ge  1$,
\begin{align*}
\int_0^{2\pi } {\vert 1 \!+ \e^{i\beta }\vert ^p \d\beta }
\int_0^{2\pi } {| {D_\alpha P(\e^{i\theta })} |} ^p \d\theta \le
2\pi n^p(\vert \alpha \vert + 1)^p\int_0^{2\pi } {| {P(\e^{i\theta
})} |^p \d\theta },
\end{align*}
which immediately leads to (15) and this completes the proof of
Theorem~2.}
\end{pot}

\begin{pot}
{\rm Since $P(z)$ is a self-inversive polynomial of degree $n$, we
have $P(z)=Q(z)$ where $Q(z) = z^n\overline {P(1/\bar {z})}$.
Therefore, for every complex number $\alpha$,
\begin{equation*}
| {D_\alpha P(z)} | = | {D_\alpha Q(z)} |\quad\hbox{for all} \ z
\in  C,
\end{equation*}
so that
\begin{align*}
|D_\alpha Q(\e^{i\theta })/D_{\alpha} P(\e^{i\theta })|= 1.
\end{align*}
Using this in (22) and proceeding similarly as in the proof of
Theorem~2, we get (16) and this proves Theorem~3.}
\end{pot}

\end{document}